\documentclass[12pt]{amsart}
  
\usepackage{amsfonts,graphics,amsmath,amsthm,amsfonts,amscd,amssymb,amsmath,latexsym,multicol}
\usepackage{epsfig,url}
\usepackage{flafter}

\makeatletter

\def\jobis#1{FF\fi
  \def\predicate{#1}%
  \edef\predicate{\expandafter\strip@prefix\meaning\predicate}%
  \edef\job{\jobname}%
  \ifx\job\predicate
}

\makeatother

\if\jobis{proposal}%
\else
\fi

 \usepackage[matrix, arrow]{xy}

\DeclareMathOperator{\Supp}{Supp}

 \newcommand{\N}{\mathbb N}
 
 \newcommand{\Q}{\mathbb Q}
 \newcommand{\R}{\mathbb R}
 \newcommand{\Z}{\mathbb Z}
 \newcommand{\bir}{\dashrightarrow}
 \newcommand{\rddown}[1]{\left\lfloor{#1}\right\rfloor} 
 \newcommand{\ep}{\varepsilon}

 \numberwithin{equation}{subsection}
 \numberwithin{footnote}{subsection}

 \newtheorem{cor}[subsection]{Corollary}
 \newtheorem{lem}[subsection]{Lemma}
 
 \newtheorem{thm}[subsection]{Theorem}

{
    \newtheoremstyle{upright}%
        {8pt plus2pt minus4pt}%
        {8pt plus2pt minus4pt}%
        {\upshape}%
        {}%
        {\bfseries\scshape}%
        {}%
        {1em}%
        {}%
\theoremstyle{upright}

 \newtheorem{defn}[subsection]{Definition}

 \newtheorem{quest}[subsection]{Question}
 \newtheorem{rem}[subsection]{Remark}

}

\title[Log minimal moldels and Zariski decompositions]{On existence of log minimal models and weak Zariski decompositions}
\author{Caucher Birkar}
\date{\today}
\begin{document}
\maketitle

\begin{abstract}
We first introduce a weak type of Zariski decomposition in higher dimensions: an $\R$-Cartier divisor 
has a weak Zariski decomposition if birationally and in a numerical sense it can be written 
as the sum of a nef and an effective $\R$-Cartier divisor. We then prove that 
there is a very basic relation between Zariski decompositions and log minimal models 
which has long been expected: we prove that 
assuming the log minimal model program in dimension $d-1$, a lc pair $(X/Z,B)$ of 
dimension $d$ has a log minimal model 
if and only if $K_X+B$ has a weak Zariski decomposition$/Z$. 
\end{abstract}


\section{Introduction}
We work over an algebraically closed field $k$ of characteristic zero, and often in the relative situation, that is, when we have a 
projective morphism $X\to Z$ of normal quasi-projective varieties written as $X/Z$. See section 2 for notation and terminology.

Zariski decomposition is a central notion in birational geometry and it is not a surprise that there are several 
definitions of it. It all started with Zariski [\ref{Zariski}] who proved that any effective divisor $D$ on a smooth projective 
surface $X$ can be decomposed as $D=P+N$ where 

\begin{enumerate}
\item $P$ and $N$ are $\Q$-divisors, 
\item $P$ is nef, $N\ge 0$, 
\item $N=0$ or the intersection matrix $[(C_i\cdot C_j)]_{1\le i\le n, 1\le j\le n}$ is negative definite where $C_1,\dots, C_n$ are the 
components of $N$, and 
\item $P\cdot C_i=0$ for any $i$.   
\end{enumerate}

Fujita [\ref{Fujita}] generalised Zariski's result to the situation where $D$ is a pseudo-effective $\R$-divisor 
in which case one needs to replace (1) with: $P$ and $N$ are $\R$-divisors.

There have been many attempts to generalise the above decomposition to higer dimensions. Here we mention
two of them which are more relevant. 

\begin{defn}[Fujita-Zariski decomposition] 
Let $D$ be an $\R$-Cartier divisor on a normal variety $X/Z$. A Fujita-Zariski decomposition$/Z$ 
for $D$ is an expression $D=P+N$ such that 
\begin{enumerate}
\item $P$ and $N$ are $\R$-Cartier divisors, 
\item $P$ is nef$/Z$, $N\ge 0$, and 
\item if $f\colon W\to X$ is a projective birational morphism from a normal variety, and $f^*D=P'+N'$ 
with $P'$ nef$/Z$ and $N'\ge 0$, then $P'\le f^*P$.
\end{enumerate}
\end{defn}

This decomposition is unique if it exists.
Fujita [\ref{Fujita2}] related this type of decomposition to the finite generation of the canonical rings of elliptic $3$-folds.

\begin{defn}[CKM-Zariski decomposition] 
Let $D$ be an $\R$-Cartier divisor on a normal variety $X/Z$. A Cutkosky-Kawamata-Moriwaki-Zariski 
(CKM-Zariski for short) decomposition$/Z$ for $D$ is an expression $D=P+N$ such that 
\begin{enumerate}
\item $P$ and $N$ are $\R$-Cartier divisors, 
\item $P$ is nef$/Z$, $N\ge 0$, and 
\item the morphisms $\pi_*\mathcal{O}_X(\rddown{mP})\to \pi_*\mathcal{O}_X(\rddown{mD})$ are isomorphisms for all 
$m\in \N$ where $\pi$ is the given morphism $X\to Z$.
\end{enumerate}
\end{defn}

If the CKM-Zariski decomposition exists, then it is unique when $D$ is big$/Z$ but not necessarily unique in other cases.
There are examples of big divisors with no Zariski decomposition by either definitions [\ref{Nakayama}, Section 5.2]. 
Prokhorov [\ref{Prokhorov}] provides a nice survey of various kinds of Zariski decompositions.\\

{\textbf{Relation with log minimal models.}} It has long been  predicted that existence of log minimal models should be closely related to Zariski decompositions (cf. [\ref{KMM}, Section 7.3]). 
Let $(X/Z,B)$ be a lc pair. It is well-known that if we have a log minimal model for this pair, then birationally  there is a Zariski 
decomposition for $K_X+B$  according to both definitions of Zariski decomposition. However, the converse is known only 
in some special cases. Based on works of Moriwaki [\ref{Moriwaki}], Kawamata [\ref{Kawamata}] proved that for a klt pair $(X/Z,B)$ with 
$K_X+B$ being $\Q$-Cartier and big$/Z$, existence of 
a CKM-Zariski decomposition for $K_X+B$ implies existence of a log canonical model for $(X/Z,B)$. Shokurov [\ref{pl-flips}, Section 3] 
tried to use similar ideas to construct log flips but it did not work and he developed different methods to deal with log flips. 
The problem with both definitions is that it is frequently very hard to construct such decompositions because of condition (3). 

In this paper, by introducing a very weak type of Zariski decomposition, 
we show that there is a much more basic relation between such decompositions and  
log minimal models.

\begin{defn}[Weak Zariski decomposition]\label{d-WZD}
Let $D$ be an $\R$-Cartier divisor on a normal variety $X/Z$. A \emph{weak Zariski decomposition$/Z$} for $D$
consists of a projective birational morphism $f\colon W\to X$ from a normal variety, and a numerical equivalence 
$f^*D\equiv P+N/Z$ such that 
\begin{enumerate}
\item $P$ and $N$ are $\R$-Cartier divisors, 
\item $P$ is nef$/Z$, and $N\ge 0$.
\end{enumerate}
\end{defn}

Note that this is much weaker than the other definitions. For example, if $D\ge 0$, then by taking $f$ to be the 
identity, $P=0$ and $N=D$ we already 
have a weak Zariski decomposition. Of course, a natural thing to ask is the following

\begin{quest}\label{question}
Does every pseudo-effective$/Z$ $\R$-Cartier divisor $D$ on a normal variety $X/Z$ have a weak Zariski decomposition$/Z$?
\end{quest}

This follows from a question posed by Nakayama  [\ref{Nakayama}, Problem, page 4], when $Z=pt$, which states that: 
  let $D$ be a pseudo-effective $\R$-Cartier divisor on a normal projective variety $X$; is there a projective birational 
morphism $f\colon W\to X$ such that $P_\sigma(f^*D)$ is nef? If the answer is yes then $f^*D=P_\sigma(f^*D)+N_\sigma(f^*D)$ 
where $N_\sigma(f^*D)\ge 0$ hence we get a weak Zariski decomposition. 

Now we come to the main result of this paper.

\begin{thm}\label{induction-II}
Assume the LMMP for $\Q$-factorial dlt pairs in dimension $d-1$. Let $(X/Z,B)$ be a lc pair of dimension $d$.
Then, the following are equivalent:
\begin{enumerate}
\item $K_X+B$ has a weak Zariski decomposition$/Z$,
\item $K_X+B$ birationally has a CKM-Zariski decomposition$/Z$,
\item  $K_X+B$ birationally has a Fujita-Zariski decomposition$/Z$,
\item $(X/Z,B)$ has a log minimal model.
\end{enumerate}
\end{thm}

Perhaps the most important aspect of Theorem \ref{induction-II} is that it offeres an inductive approach to the 
minimal model conjecture, and it unifies the recent methods of [\ref{B}][\ref{BCHM}][\ref{BP}] and [\ref{ordered}][\ref{B2}].
In [\ref{B}][\ref{B-II}], we reduced the minimal model conjecture to the nonvanishing conjecture. One of the motivations 
behind the above theorem is to reduce the minimal model conjecture to the existence of weak Zariski decompositions 
for log canonical divisors which is a much weaker problem than the nonvanishing conjecture.

The following corollary gives some degree of flexibility in construction of log minimal models.

\begin{cor}\label{c-boundaries}
Assume the LMMP for $\Q$-factorial dlt pairs in dimension $d-1$. Let $(X/Z,B)$ and $(X/Z,B')$ be lc pairs of dimension $d$ 
such that $B'\le B$. If $(X/Z,B')$ has a log minimal model, then $(X/Z,B)$ also has a log minimal model.
\end{cor} 

This corollary gives yet another proof of the following result which was first proved by Shokurov [\ref{ordered}] and a simplified proof of it  was given in [\ref{B2}] in the klt case. It is also proved in [\ref{B}] but for effective pairs. 

\begin{cor}\label{c-dim4}
Let $(X/Z,B)$ be a lc pair of dimension $4$. Then, it has a log minimal model or a Mori fibre space.
\end{cor}

Finally, we remark some other possible consequences of Corollary \ref{c-boundaries}. Assume the LMMP for $\Q$-factorial dlt 
pairs in dimension $d-1$ and let $(X/Z,B)$ be a lc pair of dimension $d$. We may assume 
that $(X/Z,B)$ is log smooth. If $K_X+B$ is not pseudo-effective$/Z$, then $(X/Z,B)$ has a Mori fibre space by 
[\ref{BCHM}] or [\ref{BP}]. Now assume that $K_X+B$ is pseudo-effective$/Z$.
Put $B'=\lambda B$ where $\lambda$ is the smallest nonnegative real number 
such that $K_X+\lambda B$ is pseudo-effective$/Z$. If $B'=0$, then Corollary \ref{c-boundaries} says that 
it is enough to construct a log minimal model of $(X/Z,0)$ in the category of terminal singularities. But 
if $B'\neq 0$, then  it seems likely that one can prove that $K_X+B'\equiv M/Z$ for 
some $M\ge 0$ assuming the bundance conjecture in lower dimensions - this was pointed out by McKernan in a seminar. 
In that case we can apply Corollary \ref{c-boundaries} again.


\section{Basics}

Let $k$ be an algebraically closed field of characteristic zero fixed throughout the paper. 

A \emph{pair} $(X/Z,B)$ consists of normal quasi-projective varieties $X,Z$ over $k$, an $\R$- divisor $B$ on $X$ with
coefficients in $[0,1]$ such that $K_X+B$ is $\mathbb{R}$-Cartier, and a projective 
morphism $X\to Z$. For a prime divisor $D$ on some birational model of $X$ with a
nonempty centre on $X$, $a(D,X,B)$
denotes the log discrepancy.

A pair $(X/Z,B)$ is called \emph{pseudo-effective} if $K_X+B$ is pseudo-effective/$Z$, that is,
up to numerical equivalence/$Z$ it is the limit of effective $\R$-divisors. The pair is called \emph{effective} if
$K_X+B$ is effective/$Z$, that is, there is an $\R$-divisor $M\ge 0$ such that $K_X+B\equiv M/Z$;

Let $(X/Z,B)$ be a lc pair. By a \emph{log flip}$/Z$ we mean the flip of a $K_X+B$-negative extremal flipping contraction$/Z$ [\ref{B}, Definition 2.3], 
and by a \emph{pl flip}$/Z$ we mean a log flip$/Z$ such that $(X/Z,B)$ is $\Q$-factorial dlt and the log flip is also an $S$-flip for 
some component $S$ of $\rddown{B}$.

A \emph{sequence of log flips$/Z$ starting with} $(X/Z,B)$ is a sequence $X_i\bir X_{i+1}/Z_i$ in which  
$X_i\to Z_i \leftarrow X_{i+1}$ is a $K_{X_i}+B_i$-flip$/Z$, $B_i$ is the birational transform 
of $B_1$ on $X_1$, and $(X_1/Z,B_1)=(X/Z,B)$.
\emph{Special termination} means termination near $\rddown{B}$ of any sequence of log flips$/Z$ 
starting with a pair $(X/Z,B)$, that is, 
the log flips do not intersect $\rddown{B}$ after finitely many of them.  

\begin{defn}
Let  $(X/Z,B)$ be a pair, $f\colon W\to X$ a projective birational morphism from a normal variety, and $N$ an 
$\R$-divisor on $W$. Define $\theta(X/Z,B,N)$ to be the number of components of $f_*N$ which are not components of 
$\rddown{B}$. This is a generalisation of [\ref{B}, Definition 2.2].
\end{defn}

\begin{defn}\label{d-model}
A pair $(Y/Z,B_Y)$ is a \emph{log birational model} of $(X/Z,B)$ if we are given a birational map
$\phi\colon X\bir Y/Z$ and $B_Y=B^\sim+E$ where $B^\sim$ is the birational transform of $B$ and 
$E$ is the reduced exceptional divisor of $\phi^{-1}$, that is, $E=\sum E_j$ where $E_j$ are the
exceptional/$X$ prime divisors on $Y$. The log birational model $(Y/Z,B_Y)$ is called a \emph{log smooth model} 
of $(X/Z,B)$ if $\phi^{-1}$ is a log resolution of $(X/Z,B)$.

A log birational model $(Y/Z,B_Y)$ is a \emph{nef model} of 
$(X/Z,B)$ if in addition\\\\
(1) $(Y/Z,B_Y)$ is $\Q$-factorial dlt, and\\
(2) $K_Y+B_Y$ is nef/$Z$.\\

And  we call a nef model $(Y/Z,B_Y)$ a \emph{log minimal model} of $(X/Z,B)$ if in addition\\\\
(3) for any prime divisor $D$ on $X$ which is exceptional/$Y$, we have
$$
a(D,X,B)<a(D,Y,B_Y)
$$
\end{defn}

\begin{defn}[Mori fibre space]
A log birational model $(Y/Z,B_Y)$ of a lc pair $(X/Z,B)$ is called a Mori fibre space if 
$(Y/Z,B_Y)$ is $\Q$-factorial dlt, there is a $K_Y+B_Y$-negative extremal contraction $Y\to T/Z$ 
with $\dim Y>\dim T$, and 
$$
a(D,X,B)\le a(D,Y,B_Y)
$$
for  any prime divisor $D$ (on birational models of $X$) and the strict inequality holds if $D$ is on $X$ and contracted$/Y$.
\end{defn}

Our definitions of log minimal models and Mori fibre spaces are slightly different 
from the traditional ones, the difference being that we do not assume that $\phi^{-1}$ does not contract divisors. 
Even though we allow $\phi^{-1}$ to have exceptional divisors but these divisors are very special; if $D$ 
is any such prime divisor, then $a(D,X,B)=a(D,Y,B_Y)=0$.
Actually, in the plt case, our definition of log minimal models and the traditional one coincide (see [\ref{B}, Remark 2.6]).

\begin{rem}\label{r-mmodels}
Let $(X/Z,B)$ be a lc pair.\\

(i) Suppose that  $(W/Z,B_W)$ is a log smooth model of $(X/Z,B)$ and $(Y/Z,B_Y)$ 
a log minimal model of 
$(W/Z,B_W)$. We can also write $(Y/Z,B_Y)$ as $(Y/Z,B^\sim+E)$ where $B^\sim$  is the birational transform of 
$B$ on $X$ and $E$ is the reduced divisor whose components are the exceptional$/X$ divisors on $Y$. 
 Let $D$ be a prime divisor on $X$ contracted$/Y$. Then, 
$$
a(D,X,B)=a(D,W,B_W)<a(D,Y,B_Y)=a(D,Y,B^\sim+E)
$$
which implies that $(Y/Z,B_Y)$ is a log minimal model of $(X/Z,B)$. 

Now assume that $f^*(K_X+B)\equiv P+N/Z$ where $P$ is nef$/Z$, $N\ge 0$, and $f$ is the given morphism $W\to X$. 
Then, there is 
an exceptional$/X$ $\R$-divisor $F\ge 0$ such that $K_W+B_W=f^*(K_X+B)+F$ and we
get a weak Zariski decomposition$/Z$ as $K_W+B_W\equiv P_W+N_W/Z$ where $P_W=P$, $N_W=N+F$ 
and we have 
$$
\theta(X/Z,B,N_W)=\theta(X/Z,B,N)=\theta(W/Z,B_W,N)=\theta(W/Z,B_W,N_W)
$$ 
 because $f_*N_W=f_*N$ and every prime exceptional$/X$ divisor on $W$ is a component of  $\rddown{B_W}$, 
in particular, every component of $F$ is a component of $\rddown{B_W}$.\\

(ii) Let $(Y/Z,B_Y)$ be a log minimal model of $(X/Z,B)$ 
and take a common resolution $f\colon W\to X$ 
and $g\colon W\to Y$. Then by applying the negativity lemma to 
$$
N:=f^*(K_X+B)-g^*(K_Y+B_Y)=\sum_D a(D,Y,B_Y)D-a(D,X,B)D
$$
with respect to $f$ we see that $N$ is effective, and exceptional$/Y$ where $D$ runs over the 
prime divisors on $W$.
This, in particular, means that $a(D,X,B)\le a(D,Y,B_Y)$ for any prime divisor $D$ on birational models of $X$.\\

(iii) If $(X/Z,B)$ is plt and $(Y/Z,B_Y=B^\sim+E)$ is a log minimal model, then $E=0$ otherwise for any component $D$ of $E$ on $Y$, 
$0<a(D,X,B)\le a(D,Y,B_Y)=0$ which is not possible.
\end{rem}

\begin{defn}\label{d-1}
For an $\R$-divisor $D=\sum d_iD_i$, let $D^{\le 1}:=\sum d_i'D_i$ where $d_i'=\min\{d_i,1\}$. As usual, $D_i$ are distinct 
prime divisors.
\end{defn}


\section{Proof of results}

The following result on extremal rays is an important ingredient of the proof of Theorem \ref{induction-II}.

\begin{lem}\label{l-ray}
Let $(X/Z,B)$ be a $\Q$-factorial dlt pair such that 
\begin{enumerate}
\item $K_X+B\equiv P+N/Z$,
\item $P$ is nef$/Z$, $N\ge 0$, and 
\item $\theta(X/Z,B,N)=0$. 
\end{enumerate}
Put
$$
\mu=\sup \{t\in [0,1] \mid P+tN ~~~\mbox{is nef$/Z$}\}. 
$$
Then, either $\mu=1$ in which case $K_X+B$ is nef$/Z$ or $\mu<1$ in which case there is a $K_X+B$-negative extremal ray $R/Z$ 
such that $N\cdot R<0$ and $(P+\mu N)\cdot R=0$.
\end{lem}
\textbf{Proof.}
If $\mu=1$ then obviously $K_X+B$ is nef$/Z$. So, assume that $\mu<1$. By construction,  $P+\mu N$ is nef$/Z$ but 
$P+(\mu+\ep')N$ is not nef$/Z$ for any $\ep'>0$. In particular, for any $\ep'>0$ there is a $K_X+B$-negative extremal ray $R/Z$ 
such that $(P+(\mu+\ep')N)\cdot R<0$ and $(P+(\mu+\ep)N)\cdot R=0$ for some $\ep\in [0,\ep')$.

If  there is no $K_X+B$-negative extremal ray $R/Z$ such that $(P+\mu N)\cdot R=0$, then there is an infinite  
strictly decreasing 
sequence of sufficiently small positive real numbers $\ep_i$ and $K_X+B$-negative extremal rays $R_i/Z$ 
such that $\lim_{i\to \infty} \ep_i=0$ and
$(P+(\mu+\ep_i)N)\cdot R_i=0$. For each $i$, let $\Gamma_i$ be an extremal curve for $R_i$ 
(cf. [\ref{ordered}, Definition 1] or [\ref{B-II}, section 3]).
By [\ref{B-II}, section 3], 
$$
-2\dim X\le (K_X+B)\cdot \Gamma_i=
(P+(\mu+\ep_i)N)\cdot \Gamma_i+(1-\mu-\ep_i)N\cdot \Gamma_i
$$
$$
=(1-\mu-\ep_i)N\cdot \Gamma_i<0
$$ 
hence $N\cdot \Gamma_i$ is bounded from above and from below. Therefore, 
$$
\lim_{i\to \infty} (P+\mu N)\cdot \Gamma_i=\lim_{i\to \infty} -\ep_iN\cdot \Gamma_i=0
$$
Since $\theta(X/Z,B,N)=0$, every component of $N$ is also a component of $B$. 
Hence there is $\delta$ such that $0<\delta \ll \mu$, $(X,B-\delta N)$ is dlt, $(K_X+B-\delta N)\cdot \Gamma_i<0$, 
and $\Supp (B-\delta N)=\Supp B$. 
We can write
$$
K_X+B-\delta N\equiv P+\mu N+(1-\mu-\delta)N/Z
$$  

Assume that
we have a component $D$ of $(1-\delta-\mu)N$ with non-rational coefficient. If $D\cdot \Gamma_i\ge 0$ for infinitely many $i$, then 
decrease the coefficient of $D$ in $(1-\delta-\mu)N$ very little so that the coefficient becomes 
rational. Otherwise, increase the coefficient of $D$ in $(1-\delta-\mu)N$ very little so that the coefficient becomes 
rational. Continuing this process, we can construct a $\Q$-divisor $N'\ge 0$ and  a dlt pair 
$(X/Z,B')$ so that $K_X+B'\equiv P+\mu N+N'/Z$ with $||B-B'||$ sufficiently small and $(K_X+B')\cdot \Gamma_i<0$ 
for any $i\in I$ where $I\subseteq \N$ is some infinite subset. Moreover, similar arguments as above 
show that $(N'-\ep_iN)\cdot \Gamma_i$ is bounded from above and from below, for any $i\in I$. 
 Since $N\cdot \Gamma_i$ is bounded, we get the  
boundedness of $N'\cdot \Gamma_i$ for $i\in I$ which implies that there are only finitely many such intersection numbers 
since $N'$ is rational. 

Now by  [\ref{B-II}, Remark 3.1](or [\ref{ordered}, Proposition 1]), there are positive real numbers $r_1,\dots, r_l$ and $m\in \N$ so that for any 
$i$ there are $n_{i,1},\dots, n_{i,l}\in \Z$ satisfying $n_{i,j}\ge -2m\dim X$ and such that 
$$
(K_X+B')\cdot \Gamma_i=\sum_j \frac{r_jn_{i,j}}{m}
$$ 
For $i\in I$, by considering the inequalities $-2\dim X\le (K_X+B')\cdot \Gamma_i<0$ 
we deduce that there can be only finitely 
many possibilities for the $n_{i,j}$ and hence for the intersection numbers $(K_X+B')\cdot \Gamma_i$.
 This is a contradiction as 
$$
0<(P+\mu N)\cdot \Gamma_i=(K_X+B')\cdot \Gamma_i-N'\cdot \Gamma_i
$$ 
and $\lim_{i\to \infty} (P+\mu N)\cdot \Gamma_i=0$.
$\Box$\\\\

\begin{defn}[LMMP using a weak Zariski decomposition]\label{d-lmmp-nef}
Let $(X/Z,B)$ be a $\Q$-factorial dlt pair such that $K_X+B\equiv P+N/Z$ where $P$ is nef$/Z$, $N\ge 0$, and 
$\theta(X/Z,B,N)=0$. If $K_X+B$ is not nef$/Z$, then as in Lemma \ref{l-ray} there exist $\mu\in [0,1)$ and an extremal ray $R/Z$ such that 
$N\cdot R<0$, $(P+\mu N)\cdot R=0$ and $P+\mu N$ is nef$/Z$, in particular, 
$(K_X+B)\cdot R<0$. By replacing $P$ with $P+\mu N$ we may assume that 
$P\cdot R=0$. Assume that $R$ defines a divisorial contraction or a log flip 
$X\bir X'/Z$, and let $B',P'$, and $N'$ be the birational transforms of $B,P$, and $N$ respectively. Obviously, $(X'/Z,B')$ is $\Q$-factorial dlt, 
$K_{X'}+B'\equiv P'+N'/Z$ 
where $P'$ is nef$/Z$, $N'\ge 0$, and $\theta(X'/Z,B',N')=0$. 
We can apply Lemma \ref{l-ray} again  and so on. In this way we obtain a particular kind of LMMP 
which we will refer to as the  \emph{LMMP using a weak Zariski decomposition} or more specifically 
the \emph{LMMP$/Z$ on $K_X+B$ using $P+N$}. 
\end{defn}

{\flushleft \textbf{Proof of Theorem \ref{induction-II}.}}
$(4) \implies (3)$: Let $(Y/Z,B_Y)$ be a log minimal model of $(X/Z,B)$. 
By definition, $B_Y=B^\sim +E$ where $B^\sim$ is the birational transform of $B$ and $E$ is the 
reduced divisor whose components are the exceptional$/X$ prime divisors of $Y$. 
Then, by Remark \ref{r-mmodels} (ii) there is a common resolution $g\colon V\to X$ and $h\colon V\to Y$ 
such that 
$$
g^*(K_X+B)=h^*(K_Y+B_Y)+N
$$ 
where $N$ is effective, and exceptional$/Y$. Moreover, by Definition \ref{d-model} (3) the support of $N$ 
contains the birational transform of the prime divisors on $X$ which 
are exceptional$/Y$.  Put $P=h^*(K_Y+B_Y)$ which is nef$/Z$. We prove that 
$P+N$ is a Fujita-Zariski decomposition$/Z$ of $g^*(K_X+B)$. Let $f\colon W\to V$ be a projective birational 
morphism from a normal variety and $f^*g^*(K_X+B)=P'+N'$ where $P'$ is nef$/Z$ and $N'\ge 0$. Now $f^*P-P'=N'-f^*N$ is 
antinef$/Y$ since $f^*P\equiv 0/Y$, and $(hf)_*(N'-f^*N)=(hf)_*N'\ge 0$ since $f^*N$ is exceptional$/Y$. 
Thus, by the negativity lemma  $f^*P-P'\ge 0$ hence the claim.\\

$(4) \implies (2)$: Suppose that $g\colon V\to X$, $h\colon V\to Y$, $P$, and $N$ are as above. 
We prove that the expression $g^*(K_X+B)=P+N$ is a CKM-Zariski decomposition$/Z$. 
Let $U\subseteq Z$ be an open subset and $U^{-1}$ its inverse 
image in $V$. Since $mg^*(K_X+B)-mP=mN\ge 0$, we have $\rddown{mg^*(K_X+B)}-\rddown{mP}\ge 0$ hence 
we have the natural injective map 
$$
H^0(U^{-1},\rddown{mP}) \to H^0(U^{-1},\rddown{mg^*(K_X+B)})
$$
and we need to prove that it is a bijection, for every $m\in \N$. Pick $\sigma \in H^0(U^{-1},\rddown{mg^*(K_X+B)})$. 
Then, 
$$
(\sigma)+mP+mN=(\sigma)+mg^*(K_X+B)\ge (\sigma)+\rddown{mg^*(K_X+B)}\ge 0
$$
holds in $U^{-1}$. Thus, $h_*(\sigma)+h_*mP\ge 0$ in $h(U^{-1})$ and since $(\sigma)+mP\equiv 0/Y$ we get $(\sigma)+mP\ge 0$ 
in $U^{-1}$ which implies that $(\sigma)+\rddown{mP}\ge 0$ in $U^{-1}$ hence $\sigma \in H^0(U^{-1},\rddown{mP})$.\\

$(2) \implies (1)$ \& $(3) \implies (1)$: These are clear from the definitions.\\

$(1) \implies (4)$:  Assume that $\mathfrak{W}$ is the set of pairs $(X/Z,B)$ such that\\ 
\begin{description}
\item[L] $(X/Z,B)$ is lc of dimension $d$,  
\item[Z] $K_X+B$ has a weak Zariski decomposition$/Z$, and 
\item[N] $(X/Z,B)$ has no log minimal model.\\
\end{description}

Clearly, it is enough to show that $\mathfrak{W}$ is empty. Assume otherwise and let $(X/Z,B)\in \mathfrak{W}$ and 
let $f\colon W\to X$, $P$ and $N$ be the data given by 
a weak Zariski decomposition$/Z$ for $K_X+B$ as in Definition \ref{d-WZD}. Assume in addition that $\theta(X/Z,B,N)$ is minimal. 
By taking an appropriate log resolution and 
using Remark \ref{r-mmodels} (i) 
we may assume that $(X/Z, \Supp B+N)$ is log smooth, $W=X$ and $f$ is the identity. 

First assume that $\theta(X/Z,B,N)=0$. Run the LMMP$/Z$ on $K_X+B$ using $P+N$ as in Definition \ref{d-lmmp-nef}. 
Obviously, $N$ negatively intersectcs each 
extremal ray in the process, and since $\Supp N\subseteq \Supp \rddown{B}$, the LMMP terminates by the special 
termination and we get a log minimal model for $(X/Z,B)$ which 
contradicts the assumption that $(X/Z,B)\in \mathfrak{W}$. Let us emphasize that this is the only place 
that we use the assumption on LMMP in dimension $d-1$.

From now on we assume that $\theta(X/Z,B,N)>0$. Define
$$
\alpha:=\min\{t>0~|~~\rddown{(B+tN)^{\le 1}}\neq \rddown{B}~\}
$$
(see Definition \ref{d-1}). In particular, $(B+\alpha N)^{\le 1}=B+C$ for some $C\ge 0$ supported in $\Supp N$, and $\alpha N=C+A$ 
where $A\ge 0$ is supported in $\Supp \rddown{B}$ and has no common components with $C$. Moreover, 
$\theta(X/Z,B,N)$ is the number of components of $C$.
The pair $(X/Z,B+C)$ is  $\Q$-factorial dlt and the expression 
$$
K_X+B+C\equiv P+N+C/Z
$$ 
is a weak Zariski decomposition$/Z$. By construction
$$
\theta(X/Z,B+C,N+C)<\theta(X/Z,B,N)
$$
and so $(X/Z,B+C)\notin \mathfrak W$ as $\theta(X/Z,B,N)$ is minimal. Therefore, since $(X/Z,B+C)$ satisfies 
conditions \textbf{L} and \textbf{Z} above, it has a log minimal model,
say $(Y/Z,(B+C)_Y)$. 

Let $g\colon V\to X$ and $h\colon V\to Y$ be a common resolution. 
By definition, $K_Y+(B+C)_Y$ is nef/$Z$. 
As in the proof of $(4) \implies (3)$ the expression 
$g^*(K_X+B+C)= P'+N'$ is a Fujita-Zariski decomposition$/Z$ where $P'=h^*(K_Y+(B+C)_Y)$ 
and $N'\ge 0$ is exceptional$/Y$. On the other hand, we have the expression 
$$
g^*(K_X+B+C)= (g^*P+(P'+N')-g^*(P+N+C))+g^*(N+C)
$$
where the first part is nef$/Z$ and the second part is effective. 
So by the definition of Fujita-Zariski decomposition we have $g^*(N+C)\ge N'$. 
Therefore, $\Supp N'\subseteq \Supp g^*(N+C)=\Supp g^*N$. Now, 

\begin{equation*}
\begin{split}
(1+\alpha)g^*(K_X+B) & \equiv g^*(K_X+B)+\alpha g^*P+ \alpha g^*N\\
& \equiv g^*(K_X+B)+\alpha g^*P+g^* C+g^*A\\
& \equiv P'+\alpha g^*P+N'+g^*A/Z
\end{split}
\end{equation*}
hence we get a weak Zariski decomposition$/Z$ as $g^*(K_X+B)\equiv P''+N''/Z$ where 
$$
P''=\frac{1}{1+\alpha}(P'+ \alpha g^*P)  \mbox{\hspace{0.5cm}  and \hspace{0.5cm}}  N''=\frac{1}{1+\alpha}(N'+g^*A)
$$
and $\Supp N''\subseteq \Supp g^* N$ hence $\Supp g_*N''\subseteq \Supp N$. Since $\theta(X/Z,B,N)$ is minimal,  
$$
\theta(X/Z,B,N)=\theta(X/Z,B,N'')
$$ 
 So, every component of $C$ is also a component 
of $g_*N''$ which in turn implies that every component of $C$ is also a component of $g_*N'$. But $N'$ is exceptional$/Y$ 
hence so is $C$ which means that $(B+C)_Y=B^\sim+C^\sim+E=B^\sim+E=B_Y$ where $\sim$ stands for birational transform 
and $E$ is the exceptional divisor of $Y\bir X$. Thus, we have $P'=h^*(K_Y+B_Y)$. Even though $K_Y+B_Y$ 
is nef$/Z$ but $(Y/Z,B_Y)$ is not necessarily a log minimal model of $(X/Z,B)$ because condition (3) of Definition 
\ref{d-model} may not be satisfied.

Let $G$ be the largest $\R$-divisor such that $G\le g^*C$ and $G\le {N}'$. By letting $\tilde{C}=g^*C-G$ and 
$\tilde{N}'= N'-G$ we get an expression 
$$
g^*(K_X+B)+\tilde{C}= P'+\tilde{N}'
$$ 
where $\tilde{C}$ and $\tilde{N}'$ are effective with no common components. 

 Assume that $\tilde{C}$ is exceptional$/X$.  Then, $g^*(K_X+B)-P'=\tilde{N}'-\tilde{C}$ is antinef$/X$ 
and by the negativity lemma $\tilde{N}'-\tilde{C}\ge 0$ which implies that $\tilde{C}=0$ since $\tilde{C}$ and $\tilde{N}'$ 
have no common components. Thus,
$$
g^*(K_X+B)-h^*(K_Y+B_Y)=\sum_D a(D,Y,B_Y)D-a(D,X,B)D=\tilde{N}'
$$
where $D$ runs over the prime divisors on $V$.
If $\Supp g_*\tilde{N}'=\Supp g_*N'$, then $\Supp \tilde{N}'$ contains the birational transform of 
all the prime exceptional$/Y$ divisors on $X$ 
hence $(Y/Z,B_Y)$ is a log minimal model of $(X/Z,B)$, a contradiction. Thus, 
$$
\Supp (g_*N'-g_*G)=\Supp g_*\tilde{N}'\subsetneq \Supp g_* N'\subseteq \Supp N
$$  
so some component of $C$ is not a component of $g_*\tilde{N}'$ because $\Supp g_*G\subseteq \Supp C$. Therefore, 
$$
\theta(X/Z,B,N)>\theta(X/Z,B,\tilde{N}')
$$ 
which implies that $(X/Z,B)\notin \mathfrak{W}$, a contradiction again.

 So, we may assume that $\tilde{C}$ is not exceptional$/X$. Let 
$\beta>0$ be the smallest number such that $\tilde{A}:=\beta g^* N-\tilde{C}$ satisfies $g_*\tilde{A}\ge 0$. Then 
there is a component of $g_*\tilde{C}$ which is not a component of $g_*\tilde{A}$. 
Now 

\begin{equation*}
\begin{split}
(1+\beta)g^*(K_X+B) & \equiv g^*(K_X+B)+\beta g^*N+\beta g^*P\\
&  \equiv g^*(K_X+B)+\tilde{C}+\tilde{A}+\beta g^*P\\
& \equiv P'+\beta g^*P+ \tilde{N}'+\tilde{A}/Z
\end{split}
\end{equation*}
where $\tilde{N}'+\tilde{A}\ge 0$ by the negativity lemma. 
Thus, we get a weak Zariski decomposition$/Z$ as $g^*(K_X+B)\equiv P'''+N'''/Z$
 where 
$$
P'''=\frac{1}{1+\beta}(P'+\beta g^*P)  \mbox{\hspace{0.5cm}  and \hspace{0.5cm}}  N'''=\frac{1}{1+\beta}(\tilde{N}'+\tilde{A})
$$
and $\Supp g_*N'''\subseteq \Supp N$.
Moreover, by construction, there is a component $D$ of $g_*\tilde{C}$ which is not a component of $g_*\tilde{A}$. 
Since $g_*\tilde{C}\le C$, $D$ is a component of $C$ hence of $N$, and since $\tilde{C}$ and $\tilde{N}'$ 
have no common components, $D$ is not a component of $g_*\tilde{N}'$. Therefore, $D$ is not a component of 
$g_*N'''=\frac{1}{1+\beta}(g_*\tilde{N}'+g_*\tilde{A})$, and   
$$
\theta(X/Z,B,N)>\theta(X/Z,B,N''')
$$ 
which gives a contradiction.
$\Box$\\

{\flushleft \textbf{Proof of Corollary \ref{c-boundaries}.}} 
Let $(Y/Z,B'_Y)$ be a log minimal model of $(X/Z,B')$ and let $g\colon W\to X$ and $h\colon W\to Y$ 
be a common resolution. We can write $g^*(K_X+B')=P+N$ where $P=h^*(K_Y+B'_Y)$ and $N\ge 0$. 
Then, 
$$g^*(K_X+B)=g^*(K_X+B')+g^*(B-B')=P+N+g^*(B-B')
$$ 
is a weak Zariski decomposition$/Z$ because $N+g^*(B-B')$ is effective. 
Now apply Theorem \ref{induction-II}.
$\Box$\\

\begin{rem}\label{r-boundaries} 
Obviously, with the same arguments as in the proof of Corollary \ref{c-boundaries} 
we could prove this stronger statement: assume the LMMP in 
dimension $d-1$ for $\Q$-factorial dlt pairs, and assume that $(X/Z,B)$ and $(X'/Z,B')$ are lc pairs of 
dimension $d$ such that there is a projective birational morphism $f\colon X'\to X$ with 
the property $K_{X'}+B'\le f^*(K_X+B)$; then $(X'/Z,B')$ having a log minimal model implies that 
$(X/Z,B)$ also has a log minimal model. This idea is used in the second proof of 
Corollary \ref{c-dim4}.
\end{rem}

{\flushleft \textbf{Proof of Corollary \ref{c-dim4} in the klt case.}} 
Let $(X/Z,B)$ be a klt pair of dimension $4$. By taking a crepant $\Q$-factorial terminal model 
we may assume that the pair is $\Q$-factorial and terminal in codimension$\ge 2$. First 
assume that $K_X+B$ is not pseudo-effective$/Z$. Existence of a Mori fibre space for 
$(X/Z,B)$ is proved in [\ref{BCHM}][\ref{BP}]. 
Now assume that 
$K_X+B$ is pseudo-effective$/Z$. Let 
$$
\mathcal{E}=\{B' \mid K_X+B' ~~\mbox{is pseudo-effective$/Z$ and $0\le B'\le B$}\}
$$
which is a compact subset of the $\R$-vetcor space $V$ generated by the components of $B$. 
If $0\in \mathcal{E}$, then $K_X$ is pseudo-effective$/Z$. In this case, $(X/Z,0)$ has 
a log minimal model because termination holds for $K_X$ by [\ref{KMM}]. We can then apply Corollary \ref{c-boundaries}. 
If $0\notin \mathcal{E}$, then there is $B'\in \mathcal{E}$ such that for any $0\le B''\lneq B'$, $K_X+B''$ is 
not pseudo-effective$/Z$: $B'$ is a point in $\mathcal{E}$ which has minimal 
distance from $0$ with respect to the standard metric on $V$. 
  
Run the LMMP$/Z$ on $K_X+B'$. No component of $B'$ is contracted by the LMMP: otherwise let $X_i\bir X_{i+1}$ 
be the sequence of log flips and divisorial contractions of this LMMP where $X=X_1$. 
Pick $j$ so that $\phi_j\colon X_j\bir X_{j+1}$ is a divisorial contraction 
which contracts a component $D_j$ of $B_j'$, the birational transform of $B'$. Now there is $a>0$ 
such that 
$$
K_{X_j}+B_j'=\phi_j^*(K_{X_{j+1}}+B_{j+1}')+aD_j
$$ 
Since $K_{X_{j+1}}+B_{j+1}'$ is pseudo-effective$/Z$, 
$K_{X_{j}}+B_{j}'-aD_j$ is pseudo-effective$/Z$ which implies that $K_X+B'-bD$ is pseudo-effective$/Z$ 
for some $b>0$ where 
$D$ is the birational transform of $D_j$, a contradiction.
On the other hand, by construction $(X/Z,B')$ is terminal in codimension $\ge 2$ and 
since no component of $B'$ is contracted, being terminal in codimension$\ge 2$ is preserved under 
the LMMP. 
The LMMP terminates by [\ref{Fujino}][\ref{mld's}] and we get a log minimal model of $(X/Z,B')$. 
Apply Corollary \ref{c-boundaries} once more to finish the proof.  
$\Box$\\

\begin{defn}
For an $\R$-divisor $D$ on a normal variety $X/Z$ define the stable base locus of $D$ over 
$Z$, denoted by ${\bf{B}}(D/Z)$, to be the complement of the set of points $x\in X$ such 
that $x\notin \Supp D'$ for some $D'\ge 0$ with the 
property $D'\sim_\R D/\pi(x)$ where $\pi$ is the given morphism $X\to Z$ and the condition 
$D'\sim_\R D/\pi(x)$ means that $D'\sim_\R D$ over some neighbourhood of $\pi(x)$.
\end{defn}

\begin{rem}\label{r-stable}
If $x\in X$ and $D=\lim_{i\to \infty} D_i$ in $N^1(X/Z)$ where 
$x\notin {\bf{B}}(D_i/Z)$ for $i\gg 0$, then there cannot be a curve $C$ contracted$/Z$ and passing through 
$x$ with $D\cdot C<0$ otherwise $D_i\cdot C<0$ for any $i\gg 0$ and $C\subseteq {\bf{B}}(D_i/Z)$ 
which contradicts the assumption $x\notin {\bf{B}}(D_i/Z)$.\\
\end{rem}

{\flushleft \textbf{Proof of Corollary \ref{c-dim4} in the lc case.}}
This (sketchy) proof follows some of the ideas in [\ref{ordered}][\ref{B2}].
Let $(X/Z,B)$ be a lc pair of dimension $4$. 
We may assume that it is $\Q$-factorial dlt by replacing it with a crepant $\Q$-factorial dlt model. 
Moreover, we may assume that any LMMP$/Z$ on $K_X+B$ will not contract divisors 
and that any step will not intersect $\Supp \rddown{B}$.
Let $H\ge 0$ be an ample$/Z$ divisor so that $K_X+B+H$ is dlt and nef$/Z$. 
Run the LMMP/$Z$ on $K_X+B$ with scaling of 
$H$. If the LMMP terminates, then we get a log minimal model or a Mori fibre space. Otherwise,  
we get an infinite sequence $X_i \bir X_{i+1}/Z_i$ of log flips$/Z$ where 
$X=X_1$ and that none of the log flips intersect $\Supp \rddown{B}$. 

Let $\lambda_1\ge \lambda_2 \ge\dots$ be the numbers determined by the LMMP with scaling (see [\ref{B-II}, Definition 2.3]). So,   
$K_{X_i}+B_i+\lambda_i H_i$ is nef$/Z$, $(K_{X_i}+B_i)\cdot R_i<0$ and  $(K_{X_i}+B_i+\lambda_i H_i)\cdot R_i=0$ where $B_i$ and $H_i$ are the  birational transforms of 
$B$ and $H$ respectively and $R_i$ is the extremal ray which defines the flipping contraction $X_i\to Z_i$.

Put $\lambda=\lim_{i\to \infty} \lambda_i$.  
If the limit is attained, that is, $\lambda=\lambda_i$ for some $i$, then the sequence terminates by [\ref{BCHM}][\ref{BP}]. So, we 
assume that the limit is not attained. By replacing 
$B_i$ with $B_i+\lambda H_i$, we can assume that $\lambda=0$, that is, $\lim_{i\to \infty} \lambda_i=0$. Since 
$K_{X_i}+B_i+\lambda_i H_i$ is nef$/Z$, it is semi-ample$/Z$ by the base point free theorem. Though $K_{X_i}+B_i+\lambda_i H_i$
 may not be klt but we can perturb it to the klt situation using the fact that $H_1$ is ample$/Z$.

Put $\Lambda_i:=B_i+\lambda_i H_i$.  By the results of [\ref{BCHM}][\ref{BP}], we can construct a crepant model  
$(Y_i/Z,\Theta_i)$ of $(X_i/Z,\Lambda_i)$ such that $(Y_i/Z,\Theta_i)$ is $\Q$-factorial, terminal in codimension$\ge 2$ 
outside $\Supp \rddown{\Theta_i}$, and such that no exceptional divisor of 
$Y_i\to X_i$ is a component of $\rddown{\Theta_i}$. We may assume that all the $Y_i$ are isomorphic in 
codimension one. By construction, $K_{Y_i}+\Theta_i$ is 
semi-ample$/Z$, and 

\begin{equation}\label{eq-stable}
{\bf{B}}(K_{Y_1}+\Theta_i^\sim/Z)\cap\Supp \rddown{\Theta_1}=\emptyset
\end{equation}
where $\Theta_i^\sim$ is the birational transform of $\Theta_i$: in fact let $y_1\in \Supp \rddown{\Theta_1}$, 
$x_1$ its image on $X_1$, $x_i$ the corresponding point on $X_i$, $y_i\in Y_i$ whose image on $X_i$ is $x_i$, 
and $z$ the image of all these points on $Z$; 
since $K_{Y_i}+\Theta_i$ is semi-ample$/Z$, there is $D_i\ge 0$ such that $y_i\notin \Supp D_i$ 
and $D_i\sim_\R K_{Y_i}+\Theta_i/z$. Let $D_1$, $D_i'$ and $D_1'$ be the birational transforms of $D_i$ 
on $Y_1$, $X_i$ and $X_1$ respectively. Assume that $y_1\in \Supp D_1$. Then, $x_1\in \Supp D_1'$ which 
implies that $x_i\in \Supp D_i'$ and $y_i\in \Supp D_i$, a contradiction.

 Let $\Delta=\lim_{i\to \infty} \Theta_i^{\sim}$ on $Y_1$. The limit is obtained component-wise. 
Thus, $K_{Y_1}+\Delta$ is a limit of movable/$Z$ divisors which in particular means that it is pseudo-effective/$Z$.
Moreover, $\rddown{\Theta_i^\sim}=\rddown{\Theta_1}=\rddown{\Delta}$ for any $i$, and $(Y_1/Z,\Delta)$ is 
terminal in codimension$\ge 2$ outside $\Supp \rddown{\Delta}$ because $\Delta\le \Theta_1$.   
Now run the LMMP$/Z$ on $K_{Y_1}+\Delta$. No divisor will be contracted because  
$K_{Y_1}+\Delta=\lim_{i\to \infty} K_{Y_1}+\Theta_i^\sim$ and each $K_{Y_1}+\Theta_i^\sim$ 
is movable$/Z$. Moreover, none of the extremal rays in the process intersect $\Supp \rddown{\Delta}$ 
because of \eqref{eq-stable} and Remark \ref{r-stable}. 
Since $(Y_1/Z,\Delta)$ is terminal in codimension $\ge 2$ outside $\Supp \rddown{\Delta}$,  
the LMMP terminates with a log minimal model $(W/Z,\Delta_W)$ by [\ref{Fujino}][\ref{mld's}]. Now apply Corollary \ref{c-boundaries} 
(see Remark \ref{r-boundaries}).
$\Box$


\vspace{2cm}

\flushleft{DPMMS}, Centre for Mathematical Sciences,\\
Cambridge University,\\
Wilberforce Road,\\
Cambridge, CB3 0WB,\\
UK\\
email: c.birkar@dpmms.cam.ac.uk

\end{document}